\documentclass[a4paper,12pt]{article}
\usepackage{amssymb}
\usepackage{latexsym}
\usepackage{amsmath,amsthm,amssymb}
\usepackage{amsfonts}
\usepackage{eufrak}

\newtheorem{thm}{Theorem}[section]

\theoremstyle{definition}

\newtheorem{rem}[thm]{Remark}

\begin{document}

\begin{center}
{\Large Cram\'er transform of Rademacher series}
\end{center}
\begin{center}
{\sc Krzysztof Zajkowski}
\footnote{The author is supported by the Polish National Science Center, Grant no. DEC-2011/01/B/ST1/03838.}\\
Institute of Mathematics, University of Bialystok \\ 
Akademicka 2, 15-267 Bialystok, Poland \\ 
kryza@math.uwb.edu.pl 
\end{center}

\begin{abstract}
A variational formula for the Cram\'er transform of series of weighted, independent symmetric Bernoulli random variables (Rademacher series) is given.
\end{abstract}



\section{Introduction}

The {\it Cram\'er transform} defines a {\it rate function} of the {\it large deviations} for  empirical means of a sequence of i.i.d. random variables (see \cite{Cramer}). The literature
concerning much more general contexts of the {\it large deviation principles} is very vast (see for instance monographs \cite{DeStro,DemZei}). A goal of this paper is only to show some 
variational formula for the Cram\'er transform of  random variables which are  series of weighted, independent symmetric Bernoulli random variables.

The  Cram\'er transform is the {\it Legendre-Fenchel transform} of the {\it cumulant generating function} of r.v. 
We will need the general notion of the Legendre-Fenchel transform  in topological spaces (see \cite{EkeTem} or \cite{BarPre}).
Let $X$ be a real locally convex Hausdorff space and $X^*$ its dual space. By $\left\langle\cdot,\cdot\right\rangle$ we denote the canonical pairing between $X$ and $X^*$. Let
$f:X\mapsto \mathbb{R}\cup\{\infty\}$ be a function nonidentically $\infty$. 
By $\mathcal{D}(f)$ we denote   the {\it effective domain} of $f$, i.e. $\mathcal{D}(f)=\{x\in X:\;f(x)<\infty\}$.
 A function $f^\ast:X^\ast\mapsto \mathbb{R}\cup\{\infty\}$  defined by
$$
f^\ast(x^*)=\sup_{x\in X}\{\left\langle x,x^*\right\rangle-f(x)\}=\sup_{x\in \mathcal{D}(f)}\{\left\langle x,x^*\right\rangle-f(x)\}\;\;\;\;\;(x^*\in X^\ast)
$$
is called   the {\it Legendre-Fenchel transform} ({\it convex conjugate}) of  $f$ and   
a function $f^{\ast\ast}:X\mapsto \mathbb{R}\cup\{\infty\}$ defined by 
$$
f^{\ast\ast}(x)=\sup_{x^\ast\in X^\ast}\{\left\langle x,x^\ast\right\rangle-f^\ast(x^\ast)\}=\sup_{x^\ast\in \mathcal{D}(f^\ast)}\{\left\langle x,x^*\right\rangle-f^\ast(x^\ast)\}\;\;\;\;\;(x\in X)
$$
is called the {\it convex biconjugate} of $f$.

The functions $f^\ast$ and $f^{\ast\ast}$ are convex and lower semicontinuous in the weak* and weak topology on $X^\ast$ and $X$, respectively. Moreover,
the  {\it biconjugate theorem} states that the function $f:X\mapsto \mathbb{R}\cup\{\infty\}$  not identically equal to $+\infty$ is convex and lower semicontinuous if and only if 
$f=f^{\ast\ast}$.

Let $I$ be a countable set and $(\epsilon_i)_{i\in I}$ be a Bernoulli sequence, i.e. a sequence of i.i.d. symmetric r.v's taking values $\pm 1$. For ${\bf t}=(t_i)_{i\in I}\in \ell^2(I)\equiv \ell^2$ the series
$$
X_{\bf t}:=\sum_{i\in I}t_i\epsilon_i
$$
converges a.s.. Notice that for ${\bf t}\in \ell^1$
$$
\vert X_{\bf t}\vert\le \sum_{i\in I}\vert t_i\vert=\Vert {\bf t}\Vert_1,
$$
i.e. $X_{\bf t}$ is a bounded r.v. and we can define its cumulant generating function  on whole $\mathbb{R}$ that is
$$
\psi_{\bf t}(s)=\ln Ee^{sX_{\bf t}}
$$
for every $s\in \mathbb{R}$. 
Because $(\epsilon_i)_{i\in I}$ is i.i.d. Bernoulli sequence then
\begin{eqnarray*}
\psi_{\bf t}(s)&=& \ln \prod_{i\in I} Ee^{st_i\epsilon_i}\\
\; &=& \ln\prod_{i\in I}\frac{e^{st_i}+e^{-st_i}}{2}=\sum_{i\in I}\ln\cosh(st_i).
\end{eqnarray*}

Observe that 
$$
\psi^\prime_{\bf t}(s)=\sum_{i\in I}t_i\tanh(st_i).
$$
We can not derive an evident form of $\psi_{\bf t}^*$ by using the {\it classical Legendre transform} because
we can not solve (inverse the derivative $\psi_{\bf t}^\prime$) the equation
\begin{equation}
\label{eqn1}
\psi_{\bf t}^\prime(s)=\alpha
\end{equation}
and find 
$$
\psi_{\bf t}^*(\alpha)=\alpha s_\alpha-\psi_{\bf t}(s_\alpha),
$$
where $s_\alpha$ is a solution of the equation (\ref{eqn1}). 

The following theorem shows some variational expression on $\psi_{\bf t}^*$.
\begin{thm}
\label{mth}
Let $(\epsilon_i)_{i\in I}$ be a Bernoulli sequence and ${\bf t}=(t_i)_{i\in I}\in \ell^1(I)$. The Cram\'er transform
of a variable $X_{{\bf t}}=\sum_{i\in I}t_i\epsilon_i$ is given by the following variational formula
$$
\psi_{\bf t}^*(\alpha)=\min_{\substack{{\bf b}\in \mathcal{D}(\psi_1^*) \\ \sum_{i\in I}t_ib_i=\alpha}}\psi_1^*({\bf b})
$$
for $\alpha\in (-\Vert {\bf t}\Vert_1,\Vert {\bf t}\Vert_1)$ and $+\infty$ otherwise, where 
$$
\psi_1^*({\bf b})=\frac{1}{2}\sum_{i\in I}\Big[\big(1+b_i\big)\ln\big(1+b_i\big)
+\big(1-b_i\big)\ln\big(1-b_i\big)\Big]
$$ 
is
the convex conjugate of a functional $\psi_1:\ell^1\mapsto \mathbb{R}$ of the form $\psi_1({\bf t})=\ln Ee^{X_{\bf t}}$
and $\mathcal{D}(\psi_1^*)\subset\ell_\infty(I)$ denotes its effective domain.
\end{thm}
\begin{rem}
Presented in the next section proof techniques are similar, but not the same, to methods used by 
Ostaszewska and Zajkowski in \cite{OstZaj,Zaj}.
\end{rem}

\section{Proof of Theorem \ref{mth}}

We begin with an observation on the absolute value of the cumulant generating function: 
$\vert \psi_{\bf t}(s)\vert\le \vert s\vert \Vert {\bf t}\Vert_1$. A parameter ${\bf t}$  may be an arbitrary element of $\ell^1$. Formally we can define a function $\psi$ of two variables:
$$
\psi(s,{\bf t})=\psi_{\bf t}(s)=\ln Ee^{sX_{\bf t}}\quad {\rm for}\quad (s,{\bf t})\in\mathbb{R}\times \ell^1.
$$
Fixing ${\bf t}$ or $s$ we write $\psi(s,{\bf t})=\psi_{\bf t}(s)$ or $\psi(s,{\bf t})=\psi_s({\bf t})$, respectively.
First we derive $\psi_s^*$ and  next we show  how $\psi_{\bf t}^*$ is expressed by $\psi_s^*$.

In a standard way one can  check the convexity of $\psi_s$ for every $s\in\mathbb{R}$. Let ${\bf t},{\bf u}\in \ell^1$ and $\lambda\in(0,1)$
then 
\begin{eqnarray*}
\psi_s(\lambda{\bf t}+(1-\lambda){\bf u}) & =& \ln Ee^{s\sum_{i\in I}(\lambda t_i+(1-\lambda) u_i)\epsilon_i}\\
\; &=& \ln E\big[\big(e^{s\sum_{i\in I}t_i\epsilon_i}\big)^\lambda\big(e^{s\sum_{i\in I}u_i\epsilon_i}\big)^{1-\lambda}\big].
\end{eqnarray*}
Using the H\"older inequality for exponents $1/\lambda$ and $1/(1-\lambda)$ we get
$$
E\big[\big(e^{s\sum_{i\in I}t_i\epsilon_i}\big)^\lambda\big(e^{s\sum_{i\in I}u_i\epsilon_i}\big)^{1-\lambda}\big]
\le \big(Ee^{s\sum_{i\in I}t_i\epsilon_i}\big)^\lambda\big(Ee^{s\sum_{i\in I}u_i\epsilon_i}\big)^{1-\lambda}
$$
and, in consequence,
\begin{eqnarray*}
\psi_s(\lambda{\bf t}+(1-\lambda){\bf u})& \le &\lambda\ln Ee^{s\sum_{i\in I}t_i\epsilon_i} + (1-\lambda)\ln Ee^{s\sum_{i\in I}u_i\epsilon_i}\\
\; & = & \lambda\psi_s({\bf t})+(1-\lambda)\psi_s({\bf u}).
\end{eqnarray*}

Because
$\psi_s:\ell^1\mapsto \mathbb{R}$ and $(\ell^1)^* \simeq \ell_\infty$ then
$$
\psi_s^*:\ell_\infty\mapsto \mathbb{R}\cup\{+\infty\}.
$$
Let ${\bf a}=(a_i)_{i\in I}\in \ell_\infty$. By the definition of the convex conjugate we have
\begin{equation}
\label{psies}
\psi_s^*({\bf a})=\sup_{{\bf t}\in \ell^1}\Big\{\left\langle {\bf t},{\bf a}\right\rangle-\sum_{i\in I}\ln\cosh(st_i)\Big\},
\end{equation}
where $\left\langle {\bf t},{\bf a}\right\rangle=\sum_{i\in I}t_ia_i$.

Note that for $s=0$ we have
\begin{equation*}
\psi_0^*({\bf a})=\left\{ \begin{array}{lll}
0 \quad & {\rm if}&\quad {\bf a}={\bf 0},\\
+\infty  &  {\rm otherwise} .  
\end{array} \right.
\end{equation*}
Assume now that  $s\neq 0$. An expression in the curly bracket of (\ref{psies}), denote it by $w$, is concave and its partial derivatives along vector of basis
$e_i=(\delta_{ij})_{j\in I}$ in $\ell^1$ ($\delta_{ij}$ is the Kronecker delta) equal
$$
\frac{\partial}{\partial t_i}w({\bf t})=\frac{\partial}{\partial t_i}\Big(\sum_{i\in I}t_ia_i-\sum_{i\in I}\ln\cosh(st_i)\Big)=a_i-s\tanh(st_i).
$$
 
The expression $w$ is a sum of functions with separated variables $(t_i)_{i\in I}$. Concavity of each of these functions implies that the gradient 
$\nabla w({\bf t})=(a_i-s\tanh(st_i))_{i\in I}$  belongs to the subgradient $\partial w({\bf t})$ since 
$$\forall_{{\bf u}\in\ell^1}\quad w({\bf t})-w({\bf u})\le \sum_{i\in I}(t_i-u_i)[a_i-s\tanh(st_i)]=\left\langle {\bf t}-{\bf u},\nabla w({\bf t})\right\rangle.
$$

The concave function $w$ attained its maximum (global) at the point ${\bf t}$ if and only if ${\bf 0}\in \partial w({\bf t})$.  It suffices that
$$
\forall_{i\in I}\quad a_i-s\tanh(st_i)=0.
$$ 
Because $arc\tanh(x)=\frac{1}{2}\ln\frac{1+x}{1-x}$ for $\vert x \vert<1$ then the partial derivatives equal zero
when
$$
t_i=\frac{1}{2s}\ln\frac{1+\frac{a_i}{s}}{1-\frac{a_i}{s}}\quad {\rm for}\quad \Big\vert \frac{a_i}{s} \Big\vert<1.
$$
Substituting the above values of $t_i$'s
into (\ref{psies}) we get
$$
\psi_s^*({\bf a})=\frac{1}{2}\sum_{i\in I}\Big[\big(1+\frac{a_i}{s}\big)\ln\big(1+\frac{a_i}{s}\big)
+\big(1-\frac{a_i}{s}\big)\ln\big(1-\frac{a_i}{s}\big)\Big]\quad {\rm for}\;\Big\vert\frac{a_i}{s} \Big\vert<1.
$$

Look a bit closely at the effective domain of $\psi_s^*$ that is at the set
$$
\mathcal{D}(\psi_s^*)=\Big\{{\bf a}\in l_\infty:\;\psi_s^*({\bf a})<\infty\;\Big\}.
$$
The function $f(x)=(1+x)\ln(1+x)+(1-x)\ln(1-x)$ is even and $f(0)=0$. Since $\lim_{|x|\to 1^-}=2\ln 2$ we can extend its domain to the interval $[-1,1]$. One can check that $(1+x)\ln(1+x)+(1-x)\ln(1-x)\ge x^2$. It follows that
$$
\sum_{i\in I}\Big[\big(1+\frac{a_i}{s}\big)\ln\big(1+\frac{a_i}{s}\big)
+\big(1-\frac{a_i}{s}\big)\ln\big(1-\frac{a_i}{s}\big)\Big]\ge \frac{1}{s^2}\sum_{i\in I} a_i^2
$$
and $|a_i|\le |s|$. Let $\overline{B}_{\infty}({\bf 0};r)$ denote of the closed ball at the center ${\bf 0}$ and radius $r$ in the space $\ell_\infty$. The properties of $f$ gives that
$$
\mathcal{D}(\psi_s^*)\subset \overline{B}_{\infty}({\bf 0};|s|)\cap \ell^2.
$$
Let us note that $\mathcal{D}(\psi_s^*)$ is a symmetric set that is ${\bf a}\in \mathcal{D}(\psi_s^*)$ if and only if $-{\bf a}\in \mathcal{D}(\psi_s^*)$. Moreover it is symmetric with respect to each coordinates $a_i$ of ${\bf a}$.
 
Return to the function $\psi_{\bf t}$. Let us observe that 
$$
\vert \psi_{\bf t}^\prime(s)\vert=\Big\vert \sum_{i\in I} t_i\tanh(st_i)\Big\vert < \Vert {\bf t}\Vert_1
$$
and $\lim_{s\to\pm\infty}\psi_{\bf t}^\prime(s)=\pm \Vert {\bf t}\Vert_1$. 
It follows $\mathcal{D}(\psi_{\bf t}^*)=\psi_{\bf t}^\prime(\mathbb{R})=(-\Vert {\bf t}\Vert_1,\Vert {\bf t}\Vert_1)$.
Because  $\psi_{\bf t}$ is convex and continuous on $\mathbb{R}$ then, by the biconjugate theorem, we get 
$$
\psi_{\bf t}(s)=\psi_{\bf t}^{\ast\ast}(s)=\sup_{\alpha\in (-\Vert {\bf t}\Vert_1,\Vert {\bf t}\Vert_1)}\big\{\alpha s-\psi_{\bf t}^*(\alpha)\big\}.
$$
On the other hand
$$
\psi_{\bf t}(s)=\psi_s({\bf t})=\sup_{{\bf a}\in \mathcal{D}(\psi_s^*)}\Big\{\left\langle {\bf t},{\bf a}\right\rangle-\frac{1}{2}\sum_{i\in I}\Big[\big(1+\frac{a_i}{s}\big)\ln\big(1+\frac{a_i}{s}\big)
+\big(1-\frac{a_i}{s}\big)\ln\big(1-\frac{a_i}{s}\big)\Big]\Big\}.
$$
If we take ${\bf a}=s{\bf b}$ then $\psi_s^*(s{\bf b})=\psi_1^*({\bf b})$ with ${\bf b}\in D(\psi_1^*)$.
It means that we can rewrite the above variational principle as follows 
\begin{equation}
\label{psite}
\psi_{\bf t}(s)=\sup_{{\bf b}\in \mathcal{D}(\psi_1^*)}\Big\{s\left\langle {\bf t},{\bf b}\right\rangle-\frac{1}{2}\sum_{i\in I}\Big[\big(1+b_i\big)\ln\big(1+b_i\big)
+\big(1-b_i\big)\ln\big(1-b_i\big)\Big]\Big\}.
\end{equation}
Take now $\alpha=\left\langle {\bf t},{\bf b}\right\rangle$. Recall that 
$$
\sup_{{\bf b}\in \overline{B}_{\infty}({\bf 0};1)}\left\langle {\bf t},{\bf b}\right\rangle=\Vert {\bf t}\Vert_1.
$$
We show that every number in $(-\Vert {\bf t}\Vert_1,\Vert {\bf t}\Vert_1)$ is taken by the inner product $\left\langle {\bf t},{\bf b}\right\rangle$ over the set $\mathcal{D}(\psi_1^*)$.
Observe that  a vector ${\bf b}=\sum_{i\in J}r(sgn\;t_i)e_i$, where $J$ is some finite subset of $I$ and $r\in [-1,1]$, belongs to $\mathcal{D}(\psi_1^*)$ (only finite number of nonzero terms). For this vector we have
$$
\left\langle {\bf t},{\bf b}\right\rangle=r\sum_{i\in J}|t_i|.
$$
It follows that the inner product $\left\langle {\bf t},{\bf b}\right\rangle$ attains  over the set $\mathcal{D}(\psi_1^*)$ any number  belonging to the interval 
$(-\Vert {\bf t}\Vert_1,\Vert {\bf t}\Vert_1)$. 

For a fixed ${\bf t}\in \ell^1$,  intersect $\mathcal{D}(\psi_1^*)\subset \ell_\infty$ with a family of hyperplains 
$$
\Big\{{\bf b}\in \ell_\infty:\;\left\langle {\bf t},{\bf b}\right\rangle=\alpha\Big\}_{\alpha\in (-\Vert {\bf t}\Vert_1,\Vert {\bf t}\Vert_1)}.
$$
Now we can divide the supremum of (\ref{psite}) into two parts and get
\begin{eqnarray}
\label{form1}
\psi_{\bf t}(s) & = & \sup_{\alpha\in (-\Vert {\bf t}\Vert_1,\Vert {\bf t}\Vert_1)} \sup_{\substack{{\bf b}\in \mathcal{D}(\psi_1^*) \\ \left\langle {\bf t},{\bf b}\right\rangle=\alpha}}
\Big\{s\left\langle {\bf t},{\bf b}\right\rangle-\psi_1^\ast({\bf b})\Big\}\nonumber\\
\; & = & \sup_{\alpha\in (-\Vert {\bf t}\Vert_1,\Vert {\bf t}\Vert_1)}\Big\{s\alpha-
\inf_{\substack{{\bf b}\in \mathcal{D}(\psi_1^*) \\ \left\langle {\bf t},{\bf b}\right\rangle=\alpha}}\psi_1^\ast({\bf b})\Big\}.
\end{eqnarray}
Define a function 
$$
\varphi_{\bf t}(\alpha)=\inf_{\substack{{\bf b}\in \mathcal{D}(\psi_1^*) \\ \left\langle {\bf t},{\bf b}\right\rangle=\alpha}}\psi_1^\ast({\bf b}).
$$

We prove that  in the above definition of function $\varphi_{\bf t}$ an infimum over the set 
$\mathcal{D}(\psi_1^*)\cap\{{\bf b}\in \ell_\infty:\;\left\langle {\bf t},{\bf b}\right\rangle=\alpha\}$ is attained and we can replace it by a minimum over this set that is we prove
\begin{equation}
\label{phi}
\varphi_{\bf t}(\alpha)=\min_{\substack{{\bf b}\in \mathcal{D}(\psi_1^*) \\ \left\langle {\bf t},{\bf b}\right\rangle=\alpha}}\psi_1^*({\bf b})
\end{equation}
for $\alpha\in (-\Vert {\bf t}\Vert_1,\Vert {\bf t}\Vert_1)$ and $+\infty$ otherwise.

By Banach-Alaoglu theorem the closed (unit) ball $\overline{B}_{\infty}({\bf 0};1)\subset \ell_\infty\simeq (\ell^1)^\ast$ is weak* compact and for each ${\bf t}$ and $\alpha\in (-\Vert {\bf t}\Vert_1,\Vert {\bf t}\Vert_1)$ 
the hyperplain
$H_{{\bf t},\alpha}=\{{\bf b}\in \ell_\infty:\;\left\langle {\bf t},{\bf b}\right\rangle=\alpha\}$ is closed in this topology. We have that an intersection $\overline{B}_{\infty}({\bf 0};1)\cap H_{{\bf t},\alpha}$ is weak* compact. Let $\ell_0$ be the space of sequences with finite support. Obviously $\ell_0\cap \overline{B}_{\infty}({\bf 0};1)\subset \mathcal{D}(\psi_1^*)$ and 
$H_{{\bf t},\alpha}\cap \ell_0\neq\emptyset$.
We have 
$$
\forall_{{\bf t}\in \ell^1}\forall_{\alpha\in (-\Vert {\bf t}\Vert_1,\Vert {\bf t}\Vert_1)}\quad
\mathcal{D}(\psi_1^*)\cap H_{{\bf t},\alpha}\supset \overline{B}_{\infty}({\bf 0};1)\cap H_{{\bf t},\alpha}\cap \ell_0\neq\emptyset.
$$
Recall that the function $\psi_1^*$ is nonegative and lower semicontinuous in the weak* topology. By Weierstrass Theorem $\psi_1^*$ attains its  minimum in the compact set $\overline{B}_{\infty}({\bf 0};1)\cap H_{{\bf t},\alpha}$. Because an intersection of this set with the effective domain of $\psi_1^*$ is nonempty then it means that a nonegative infimum is attained
at some element in $\mathcal{D}(\psi_1^*)$. It follows that in the definition of $\varphi_{\bf t}$ we can replace the infimum
by minimum and the formula (\ref{phi}) holds.

The formula (\ref{form1}) means that $\psi_{\bf t}$ is the convex conjugate of $\varphi_{\bf t}$. To prove an equality
$\varphi_{\bf t}=\psi_{\bf t}^*$ we should show that $\varphi_{\bf t}$ is convex and lower semicontinuous.

First we check the convexity of $\varphi_{\bf t}$. Take $\alpha_1,\;\alpha_2\in (-\Vert {\bf t}\Vert_1,\Vert {\bf t}\Vert_1)$. If $\alpha_1$ or $\alpha_2$ do not belong to the interval $(-\Vert {\bf t}\Vert_1,\Vert {\bf t}\Vert_1)$ then
the value of $\varphi_{\bf t}$ at such $\alpha_k$ equals $\infty$ and the condition of convexity is trivially satisfied. Let ${\bf b}_k$ $(k=1,2)$ be vectors in $\mathcal{D}(\psi_1^*)\cap H_{{\bf t},\alpha_k}$ such that
$$
\varphi_{\bf t}(\alpha_k)=\min_{\substack{{\bf b}\in \mathcal{D}(\psi_1^*) \\ \left\langle {\bf t},{\bf b}\right\rangle=\alpha_k}}\psi_1^*({\bf b})
=\psi_1^*({\bf b}_k).
$$
Observe that for $\lambda\in (0,1)$
$$
\left\langle {\bf t},\lambda {\bf b}_1+(1-\lambda){\bf b}_2\right\rangle=\lambda \left\langle {\bf t},{\bf b}_1\right\rangle+
(1-\lambda)\left\langle {\bf t},{\bf b}_2\right\rangle
=\lambda\alpha_1+(1-\lambda)\alpha_2,
$$
that is $\lambda {\bf b}_1+(1-\lambda){\bf b}_2\in H_{{\bf t},\lambda\alpha_1+(1-\lambda)\alpha_2}$. The above and convexity of $\psi_1^*$ gives
\begin{eqnarray*}
\varphi_{{\bf t}}(\lambda\alpha_1+(1-\lambda)\alpha_2) & \le & \psi_1^*(\lambda {\bf b}_1+(1-\lambda){\bf b}_2)\\
\; & \le & \lambda\psi_1^*( {\bf b}_1)+(1-\lambda)\psi_1^*({\bf b}_2)
=\lambda\varphi_{{\bf t}}(\alpha_1)+(1-\lambda)\varphi_{{\bf t}}(\alpha_2).
\end{eqnarray*}

Now we prove the lower semicontinuity of $\varphi_{\bf t}$. Recall that $\psi_1^*$ is convex and lower semicontinuous in the weak* topology on $\ell_\infty$. It means that for any $c\in \mathbb{R}$ the set
\begin{equation}
\label{set1}
\{{\bf b}\in \ell_\infty:\;\psi_1^*({\bf b})\le c\}
\end{equation}
is weak* closed. Since $\psi_1^*\ge 0$ we can assume that $c\ge 0$. Because the above set  is contained in weak* compact unit ball $\overline{B}_{\infty}({\bf 0};1)\supset \mathcal{D}(\psi_1^*)$ then it is also compact in this topology. Consider a range of the set (\ref{set1}) by the functional $l_{{\bf t}}:=\left\langle {\bf t},\cdot\right\rangle$, i.e.
\begin{equation}
\label{set2}
l_{\bf t}\Big(\big\{\psi_1^*({\bf b})\le c\big\}\Big).
\end{equation}
Since for each ${\bf t}\in \ell^1$ the linear functional $l_{\bf t}$ is continuous  on $\ell_\infty$ (also in the weak* topology), by the intermediate and extreme value theorems we get that the set (\ref{set2}) is a closed interval.
By symmetry of the set (\ref{set1}) and linearity of the functional $l_{\bf t}$ we get the existence of a real number
$\alpha$ such that
$$
l_{\bf t}\Big(\big\{\psi_1^*({\bf b})\le c\big\}\Big)=[-\alpha,\alpha].
$$

We show that
$$
\varphi_{\bf t}^{-1}((-\infty,c])=[-\alpha,\alpha].
$$
Let  $\beta\in \varphi_{\bf t}^{-1}((-\infty,c])$. Since $\psi_1^*$ is lower semicontinuous, there exists  
${\bf b}_\beta$ such that
$$
c\ge \varphi_{\bf t}(\beta)= \min_{\substack{{\bf b}\in \mathcal{D}(\psi_1^*) \\ \left\langle {\bf t},{\bf b}\right\rangle=\beta}}\psi_1^*({\bf b})= 
\psi_1^*({\bf b}_\beta).
$$
That is $\left\langle {\bf t},{\bf b}_\beta\right\rangle=\beta\in[-\alpha,\alpha]$. 
Conversely, let $\beta\in[-\alpha,\alpha]$. Since $l_{{\bf t}}=\left\langle {\bf t},\cdot\right\rangle$ is continuous on the connected set $\{\psi_1^*({\bf b})\le c\}$,
there is ${\bf b}_\beta^\prime\in\{\psi_1^*({\bf b})\le c\}$ such that
$$
\left\langle {\bf t},{\bf b}_\beta^\prime\right\rangle=\beta.
$$
Note that
$$
\varphi_{\bf t}(\beta)= \min_{\substack{{\bf b}\in \mathcal{D}(\psi_1^*) \\ \left\langle {\bf t},{\bf b}\right\rangle=\beta}}\psi_1^*({\bf b})\le 
\psi_1^*({\bf b}_\beta^\prime)\le c,
$$
that is $\beta\in \varphi_{\bf t}^{-1}((-\infty,c])$.  

Because $\varphi_{\bf t}$ is convex and lower semicontinuos then $\psi_{\bf t}^\ast=\varphi_{\bf t}$, which completes the proof.

\begin{rem}
The result of Theorem \ref{mth} is similar to those obtained by the {\it contraction principle} (see for instance 
\cite{DemZei}) but let us emphasize that we used the space of parameters $\ell^1$ to generate the convex conjugate
of the investigated function and we did not consider any probability distribution on it.
\end{rem} 
\begin{rem}
Let us stress that  the proof of Theorem \ref{mth} contains some scheme which allow us to generate, under some assumptions of course, variational formulas on the Cram\'er transform for another series of random variables.
\end{rem}


\begin{thebibliography}{                    }

\bibitem{BarPre}
V. Barbu, T. Precupanu, {\it Convexity and Optimization in Banach Spaces}, 4th ed., Springer Monographs in Mathematics, Springer, Dordrecht, 2012.

\bibitem{Cramer}
H. Cram\'er, {\it Sur un nouveau th\'eor\`eme-limite de la th\'eorie des probabilit\'es},
Actualit\'es Scientifiques et Industrielles 736 (1938), 5-23. Colloque
consacr\'e \`a la th\'eorie des probabilit\'es, Vol. 3, Hermann, Paris.

\bibitem{DemZei}
A. Dembo, O. Zeitouni, {\it Large Deviations Techniques and Applications}. Corrected reprints of the second (1998) edition, Stochastic Modeling and Applied Probability, 38,
Springer-Verlag, Berlin,  2010.

\bibitem{DeStro}
J. D. Deuschel, D. W. Stroock. {\it Large Deviations}. Pure and Applied Mathematics, 137, Academic Press, Inc., Boston, 1989.


\bibitem{EkeTem}
I. Ekeland, R. T\'emam, {\it Convex Analysis and Variational Problems}, Translated from French. Corrected reprint of the 1976 English edition.
Classics in Applied Mathematics, 28, Society for Industrial and Applied Mathematics (SIAM), Philadelphia, PA, 1999.

 


\bibitem{OstZaj}
U. Ostaszewska, K. Zajkowski, {\it Cram\'er transform and {\bf t}-entropy}, Positivity 18 (2014), no. 2, 347-358.

\bibitem{Zaj}
K. Zajkowski, {\it Convex conjugates of analytic functions of logarithmically convex functional}, J. Convex Anal. 20 (2013), no. 1, 243-252. 

\end{thebibliography}
\end{document}